\documentclass[oneside,english,reqno]{amsart}
\usepackage{lmodern}
\usepackage[T1]{fontenc}
\usepackage[latin9]{inputenc}
\usepackage{geometry}
\usepackage{color}
\usepackage{prettyref}
\usepackage{units}
\usepackage{amsthm}
\usepackage{amsbsy}
\usepackage{amssymb}
\usepackage{mathdots}

\makeatletter

\newcommand{\noun}[1]{\textsc{#1}}

\numberwithin{equation}{section}
\numberwithin{figure}{section}
\numberwithin{table}{section}
\usepackage{enumitem}		
  \theoremstyle{remark}
  \newtheorem*{acknowledgement*}{\protect\acknowledgementname}
\theoremstyle{plain}
\newtheorem{thm}{\protect\theoremname}
  \theoremstyle{definition}
  \newtheorem{defn}[thm]{\protect\definitionname}

\usepackage{graphicx}
\usepackage{setspace}
\usepackage[hyphens]{url}
\usepackage[perpage,symbol]{footmisc}
\usepackage{multicol}
\usepackage{etoolbox}
\usepackage{bbm}
\usepackage{tensor}
\usepackage{enumitem}

\DefineFNsymbols*{lamportnostar}[math]{{(\dagger)}{(\ddagger)}{(\S)}{(\P)}{(\|)}{(\dagger\dagger)}{(\ddagger\ddagger)}}
\setfnsymbol{lamportnostar}

\newrefformat{prop}{Proposition~\ref{#1}}
\newrefformat{lem}{Lemma~\ref{#1}}
\newrefformat{sec}{§\ref{#1}}
\newrefformat{sub}{§\ref{#1}}
\theoremstyle{remark}
\newtheorem*{rems*}{Remarks}

\setitemize[0]{leftmargin=10pt,itemindent=0pt}
\setlist[enumerate]{leftmargin=*,widest=0}

\allowdisplaybreaks[2]

\renewenvironment{thebibliography}[1]{%
    \begin{oldthebibliography}{#1}%
      \setlength{\parskip}{0ex}%
      \setlength{\itemsep}{0ex}%
      \begin{small}
  }%
  {%
    \end{small}
    \end{oldthebibliography}%
  }

\makeatother

\usepackage{babel}
\usepackage{listings}
  \providecommand{\acknowledgementname}{Acknowledgement}
  \providecommand{\definitionname}{Definition}
\providecommand{\theoremname}{Theorem}

\begin{document}

\title[Hypermatrix Spectrum]{Approximating the spectrum of matrices and hypermatrices }

\author{Edinah K. Gnang }
\begin{abstract}
We describe a new method for computing generators of elimination ideals
associated with matrix and hypermatrix spectral constraints. We proceed
to derive from these generators iterative procedures for approximating
the spectral decomposition of matrices and hypermatrices.
\end{abstract}

\date{\today}

\maketitle
\tableofcontents{}

\section{\label{sec:Introduction}Introduction}

Many combinatorial optimization problems including instances of subgraph-isomorphism.
These combinatorial problems are known to be NP-hard, and often one
seeks to identify special families of graphs for which efficient combinatorial
algorithms can be devised. The complexity of ensuing algorithms is
often determined by combinatorial and algebraic properties of graph
encodings. A good illustration of this fact is provided by the graph
property of being an expander graph. The complexity of many combinatorial
algorithms for such a graph can be expressed in terms of its expansion
parameter \cite{AC88,Alo86,alon2008elementary,AM85,BL06,Che70,Chu97,li2001negative}.
As is well known, many graph properties, including the property of
being an expander graph, are closely tied to the spectral decomposition
of matrices deduced from incidence structures in the associated graph
\cite{Alo86,Che70,Chu97,spielman2011spectral}. The well known graph
adjacency matrix is constructed such that the $\left(i_{1},i_{2}\right)$
matrix entry equals $1$ if the associated graph admits a directed
edge connecting vertex $i_{1}$ to vertex $i_{2}$, and equals $0$
otherwise. Following a construction proposed in \cite{FW95} by Friedman
and Widgerson, one also associates with a graph a $\left(k-1\right)$-\emph{path
adjacency hypermatrix}. The $\left(k-1\right)$-path adjacency hypermatrix
corresponds to a $k$-th order hypermatrix (note that ordinary matrices
are second order hypermatrices) whose $\left(i_{1},i_{2}\cdots,i_{k}\right)$
entry equals $1$ if the ordered tuple $\left(i_{1},i_{2}\cdots,i_{k}\right)$
denotes a directed path of length $k-1$ in the associated graph,
and equals $0$ otherwise. Consequently, graph algebraic invariants
derived from the corresponding adjacency matrix are extended to include
polynomial relations between entries of path adjacency hypermatrices.
The polynomial relations which are retained as algebraic and combinatorial
invariants are usually the ones which are invariant under the natural
action of the permutation group. Such algebraic relations are to be
thought of as generalizations of the classical Cayley-Hamilton matrix
polynomial. Just as it is done for matrices, the algebraic varieties
defined by the algebraic and combinatorial invariants will be referred
to as the \emph{spectrum} or the \emph{spectral decomposition}. The
invariants associated with path adjacency hypermatrices enable us
to distinguish some non-isomorphic graphs with isospectral adjacency
matrices.

The spectral analysis of hypermatrices is considerably more difficult
to define \cite{chung1993laplacian,Lubotzky2013,parzanchevski2012isoperimetric,Fox2011proceedings,GER,2013arXiv1301.4590C,2013arXiv1309.2163Q,2014arXiv1405.7257L,FW95,2012arXiv1201.3424Q,GKZ}
when compared with the spectral analysis of matrices. Mesner and Bhattacharya
in \cite{MB90,MB94} introduced an $m$ operands product for $m$-th
order hypermatrices. E. Gnang, V. Retakh and A. Elgammal proposed
in \cite{GER} a generalization to hypermatrices of the notions of
Hermicity and unitarity. E. Gnang, V. Retakh and A. Elgammal also
proved in \cite{GER} a conjecture of Bhattacharya by using these
new definitions to extend the spectral decomposition to hypermatrices
of arbitrary order. In the present work we show that the spectral
decomposition introduced in \cite{GER} for a hypermatrix is mostly
determined by the spectral decomposition its minors in a similar spirit
to Cauchy's interlacing theorem. We show in the present work that
the majority of hypermatrix product proposed in the literature including
the Segre outer product, the contraction product, multilinear matrix
multiplication\cite{Lim2013}, as well as the Kerner product \cite{RK}
are either special cases of the general BM product or special cases
of the dual product to the general BM product. We present new algorithms
for deriving generators of the elimination ideals associated with
matrix and hypermatrix spectral decomposition constraints. The proposed
algorithms are based on generalization of Parsevals' identities also
known as resolution of the identity. We derive from the generators
the spectral elimination ideal iterative procedures for approximating
the spectral decomposition of matrices and hypermatrices. We extend
to even order hypermatrices the self-adjoint argument for establishing
the existence of real solutions to spectral constraints. Finally we
deduce from the the spectral decomposition of hypermatrices upper
and lower bounds for hypermatrix eigenvalues introduced by L.H. Lim
and L.Q. Qi in \cite{Lim,Qi:2005:ERS:1740736.1740799}.
\begin{acknowledgement*}
We would like to thank Andrei Gabrielov for providing guidance and
inspiration while preparing this manuscript. We would like to thank
Vladimir Retakh and Ahmed Elgammal for patiently introducing us to
the theory of hypermatrices. We are grateful to Avi Widgerson, Doron
Zeilberger, Yuval Filmus and Ori Parzanchevski for helpful discussions
and suggestions. The author was supported by the the National Science
Foundation, and is grateful for the hospitality of the Institute for
Advanced Study.
\end{acknowledgement*}

\section{\label{sec:Notation}Notation}

We describe for convenience of the reader the notation used throughout
the paper. The Hadamard product of matrices $\mathbf{A}$, $\mathbf{B}\in\mathbb{C}^{m\times n}$
noted $\mathbf{A}\circ\mathbf{B}$, yields a matrix of the same dimensions
whose entries are the product of corresponding entries of $\mathbf{A}$
and $\mathbf{B}$, 
\[
\left(i,j\right)-\mbox{th entry of }\left(\mathbf{A}\circ\mathbf{B}\right)\mbox{ is }a_{i,j}\,b_{i,j}.
\]
The vector product of $\mathbf{a}$, $\mathbf{b}\in\mathbb{C}^{n\times1}$
with the background $n\times n$ matrix $\mathbf{M}$ refers to the
bilinear form associated with $\mathbf{M}$ expressed as 
\[
\left\langle \mathbf{a},\mathbf{b}\right\rangle _{\mathbf{M}}:=\sum_{0\le k_{0},k_{1}<n}a_{k_{0}}m_{k_{0},k_{1}}b_{k_{1}},
\]
in particular it follows that 
\[
\left\langle \mathbf{a},\mathbf{b}\right\rangle :=\left\langle \mathbf{a},\mathbf{b}\right\rangle _{\mathbf{I}_{n}}=\sum_{0\le k<n}a_{k}\,b_{k}
\]
where the entries of $\mathbf{I}_{n}$ are given by 
\[
\left[\mathbf{I}_{n}\right]_{i,j}\,:=\left(\delta_{i,j}=\begin{cases}
\begin{array}{cc}
1 & \mbox{if }0\le i=j<n\\
0 & \mbox{otherwise}
\end{array}\end{cases}\right).
\]
The inner-product of $\mathbf{a}$, $\mathbf{b}\in\mathbb{C}^{n\times1}$
is $\left\langle \mathbf{a},\overline{\mathbf{b}}\right\rangle $.
For $\left\{ \mathbf{v}_{j}\right\} _{0\le j<n}\subset\mathbb{C}^{n\times1}$
we define the correlation product noted $\left\langle \mathbf{v}_{0},\mathbf{v}_{1},\cdots,\mathbf{v}_{n-1}\right\rangle $
to be 
\[
\left\langle \mathbf{v}_{0},\mathbf{v}_{1},\cdots,\mathbf{v}_{n-1}\right\rangle :=\left\langle \mathbf{1}_{n\times1},\,\mathbf{v}_{0}\circ\cdots\circ\mathbf{v}_{n-1}\right\rangle .
\]
It shall also be convenient to adopt the notation convention 
\[
\mathbf{a}^{\circ^{\alpha}}:=\left(\left(a_{k}\right)^{\alpha}\right)_{0\le k<n}.
\]
The $n$-dimensional vector $\mathbf{w}$ denotes the vector collecting
powers of the primitive $n$-th roots of unity with entries given
by 
\[
\mathbf{w}\,:=\left(w_{j}=\exp\left(\nicefrac{2\pi i\,j}{n}\right)\right)_{0\le j<n}.
\]
Finally, we associate with an arbitrary $\mathbf{v}\in\mathbb{C}^{n\times1}$
the $n\times n$ Vandermonde matrix 
\[
\left[\mbox{Vandermonde}\left(\mathbf{v}\right)\right]_{i,j}=\left(v_{j}\right)^{i}.
\]

\section{\label{sec:Complexes-and-buildings}Overview of the Bhattacharya-Mesner
algebra and its dual}

We recall here for convenience of the reader the basic elements of
the Bhattacharya-Mesner (BM) algebra proposed in \cite{MB90,MB94}
as a generalization of the algebra of matrices.
\begin{defn}
The Bhattacharya-Mesner \cite{MB90,MB94} algebra generalizes the
classical matrix product 
\[
\mathbf{B}=\mathbf{A}^{(1)}\cdot\mathbf{A}^{(2)}
\]
where $\mathbf{A}^{(1)}$, $\mathbf{A}^{(2)}$, $\mathbf{B}$ are
matrices of sizes $n_{1}\times k$, $k\times n_{2}$, $n_{1}\times n_{2}$,
respectively, 
\[
b_{i_{1},i_{2}}=\sum_{1\le\textcolor{red}{j}\le k}a_{i_{1},\textcolor{red}{j}}^{(1)}\,a_{\textcolor{red}{j},i_{2}}^{(2)},
\]
to an $m$-operand hypermatrix product noted
\[
\mathbf{B}=\mbox{Prod}\left(\mathbf{A}^{(1)},\,\cdots,\mathbf{A}^{(m)}\right),
\]
where $\mathbf{B}$ is an $n_{1}\times\cdots\times n_{m}$ hypermatrix,
for $i=1,\cdots,\left(m-1\right)$, $\mathbf{A}^{(i)}$ is a hypermatrix
whose size is obtained by replacing $n_{i+1}$ by $k$ in the dimensions
of the hypermatrix $\mathbf{B}$, and $\mathbf{A}^{(m)}$ is a $k\times n_{2}\times\cdots\times n_{m}$
hypermatrix, and 
\[
b_{i_{1},\cdots,i_{m}}=\sum_{1\le\textcolor{red}{j}\le k}a_{i_{1},\textcolor{red}{j},i_{3},\cdots,i_{m}}^{(1)}\cdots\,a_{i_{1},\cdots,i_{t},\textcolor{red}{j},i_{t+2},\cdots,i_{m}}^{(t)}\cdots\,a_{\textcolor{red}{j},i_{2},\cdots,i_{m}}^{(m)}.
\]
In the particular case of third order hypermatrix product noted 
\[
\mathbf{B}=\mbox{Prod}\left(\mathbf{A}^{(1)},\mathbf{A}^{(2)},\mathbf{A}^{(3)}\right)
\]
where $\mathbf{A}^{(1)}$, $\mathbf{A}^{(2)}$, $\mathbf{A}^{(3)}$
and $\mathbf{B}$ are hypermatrices of sizes $n_{1}\times k\times n_{3}$,
$n_{1}\times n_{2}\times k$, $k\times n_{2}\times n_{3}$ and $n_{1}\times n_{2}\times n_{3}$
respectively, 
\[
b_{i_{1},i_{2},i_{3}}=\sum_{1\le\textcolor{red}{j}\le k}a_{i_{1},\textcolor{red}{j},i_{2}}^{(1)}\,a_{i_{1},i_{2},\textcolor{red}{j}}^{(2)}\,a_{\textcolor{red}{j},i_{1},i_{2}}^{(3)}.
\]
A slight variation of the BM product was introduced in \cite{GER}.
The proposed variation of the BM product is called the general BM
product and noted 
\[
\mathbf{C}=\mbox{Prod}_{\mathbf{B}}\left(\mathbf{A}^{(1)},\,\cdots,\mathbf{A}^{(m)}\right).
\]
The resulting hypermatrix $\mathbf{C}$ is an $n_{1}\times\cdots\times n_{m}$
hypermatrix, while the dimensions of the hypermatrix $\mathbf{A}^{(i)}$
for $i=1,\cdots,\:m-1$ is obtained by replacing $n_{i+1}$ by $k$
in the dimensions of $\mathbf{C}$ and $\mathbf{A}^{(m)}$ is a hypermatrix
of size $k\times n_{2}\times\cdots\times n_{m}$ similarly to the
BM product. Crucially, the general BM product differs from the BM
product in the fact that the general BM product involves an additional
input hypermatrix. The additional input hypermatrix $\mathbf{B}$
is called the background hypermatrix and as such $\mathbf{B}$ must
be a cubic $m$-th order hypermatrix having all of it's sides of length
$k$ i.e. $\mathbf{B}$ is of dimension $k\times k\times\cdots\times k$
, 
\[
c_{i_{1},\cdots,i_{m}}=\sum_{1\le\textcolor{red}{j_{1}},\textcolor{red}{\cdots},\textcolor{red}{j_{m}}\le k}a_{i_{1},\textcolor{red}{j_{2}},i_{3},\cdots,i_{m}}^{(1)}\cdots a_{i_{1},\cdots,i_{t},\textcolor{red}{j_{t+1}},i_{t+2},\cdots,i_{m}}^{(t)}\cdots a_{\textcolor{red}{j_{1}},i_{2},\cdots,i_{m}}^{(m)}\,b_{\textcolor{red}{j_{1}},\textcolor{red}{\cdots},\textcolor{red}{j_{m}}}.
\]
Note that the original BM product is recovered by setting $\mathbf{B}$
to the Kronecker delta hypermatrix (i.e. the hypermatrix whose nonzero
entries all equal one and are located at the entries whose indices
all have the same value, in particular Kronecker delta matrices are
identity matrices).\\
Consider the product 
\[
\mathbf{D}=\mathbf{C}\circ\mbox{Prod}_{\mathbf{B}}\left(\mathbf{A}^{(1)},\,\cdots,\mathbf{A}^{(m)}\right),
\]
hence
\[
d_{i_{1},\cdots,i_{m}}=
\]
\[
c_{i_{1},\cdots,i_{m}}\sum_{1\le\textcolor{red}{j_{1}},\textcolor{red}{\cdots},\textcolor{red}{j_{m}}\le k}a_{i_{1},\textcolor{red}{j_{2}},i_{3},\cdots,i_{m}}^{(1)}\,\cdots\,a_{i_{1},\cdots,i_{t},\textcolor{red}{j_{t+1}},i_{t+2},\cdots,i_{m}}^{(t)}\,\cdots\,a_{\textcolor{red}{j_{1}},i_{2},\cdots,i_{m}}^{(m)}\,b_{\textcolor{red}{j_{1}},\textcolor{red}{\cdots},\textcolor{red}{j_{m}}}.
\]
We define the dual to the general BM product to be expressed as 
\[
d_{\textcolor{red}{j_{1}},\textcolor{red}{\cdots},\textcolor{red}{j_{m}}}=
\]
\begin{equation}
b_{\textcolor{red}{j_{1}},\textcolor{red}{\cdots},\textcolor{red}{j_{m}}}\sum_{1\le i_{1},\cdots,i_{m}\le k}a_{i_{1},\textcolor{red}{j_{2}},i_{3},\cdots,i_{m}}^{(1)}\cdots a_{i_{1},\cdots,i_{t},\textcolor{red}{j_{t+1}},i_{t+2},\cdots,i_{m}}^{(t)}\cdots a_{\textcolor{red}{j_{1}},i_{2},\cdots,i_{m}}^{(m)}\,c_{i_{1},\cdots,i_{m}}.
\end{equation}
The duality here arises from interchanging the indices in the summands
with the indices of the hypermatrix $\mathbf{C}$ in the product.
Note that in the case of matrices, the product is self dual. The product
dual to the general BM product was independently proposed by Kerner
in \cite{RK}. We also note that most hypermatrix product in the literature
including the Segre outer product, the contraction product, the multilinear
matrix multiplication\cite{Lim2013} all correspond to special instances
of the general BM product with additional constraints imposed on the
input hypermatrices.
\end{defn}

\section{\label{sec: Approximating the Spectrum} Spectral decomposition from
the spectrum of minors.}

We discuss here a general argument for reducing the spectral decomposition
of an arbitrary cubic $k$-th order hypermatrix to the decomposition
of cubic minors of same order having sides of length $k$. The hypermatrices
considered here are associated with directed and weighted $k$-uniform
hypergraph having no degenerate edges. In other words, the collection
of $k$ vertices which make up any hyperedge of the hypergraph must
be distinct. Such hypermatrices arise as $\left(k-1\right)$-path
adjacency hypermatrices of rooted trees whose edges are directed towards
the leaf nodes and away from the root. Such hypermatrices also arise
as $\left(k-1\right)$-path adjacency hypermatrices of directed acyclic
graphs.
\begin{thm}
\label{thm:Spectral decomposition} Let $H$ denote a directed weighted
$k$-uniform hypergraph having no degenerate hyperedges. Then the
spectral decomposition of its $k$-vertex sub-hypergraphs determine
the spectral decomposition of a larger $k$-uniform hypergraph on
$n{n \choose k}$ vertices which admit $H$ as a sub-hypergraph.
\end{thm}
We will first present the detail proof in the case of graphs and subsequently
extend the arguments to hypergraphs.
\begin{proof}
Theorem \ref{thm:Spectral decomposition} asserts that for a directed
and weighted graph having no loop edges noted 
\[
G_{1}\::=\left(V_{1}=\left\{ 0,1,\cdots,n-1\right\} ,E_{1}\subset V_{1}\times V_{1}\right),
\]
the spectral decomposition of its two vertex subgraph adjacency matrices
determine the spectral decomposition of the adjacency matrix of a
larger graph 
\[
G_{2}\::=\left(V_{2}=\left\{ 0,1,\cdots,n{n \choose 2}-1\right\} ,E_{2}\subset V_{2}\times V_{2}\right)
\]
which admit $G_{1}$ as a subgraph. Let $\mathbf{A}$ denote the $n\times n$
adjacency matrix of the graph $G_{2}$. We seek to determine the spectral
decomposition of the adjacency matrix of $G_{2}$ expressed as 
\[
\left(\begin{array}{cc}
\mathbf{A} & \mathbf{B}_{01}\\
\mathbf{B}_{10} & \mathbf{B}_{11}
\end{array}\right)=\left(\mathbf{U}\cdot\mbox{diag}\left\{ \boldsymbol{\mu}\right\} \right)\cdot\left(\mathbf{V}\cdot\mbox{diag}\left\{ \boldsymbol{\nu}\right\} \right)^{\dagger},\quad\mathbf{U}\cdot\mathbf{V}^{\dagger}=\mathbf{I}_{n{n \choose 2}}
\]
where the sub matrices $\mathbf{B}_{01}$, $\mathbf{B}_{10}$ and
$\mathbf{B}_{11}$ are matrices of size respectively given by $n\times\left(n-1\right){n \choose 2}$,
$\left(n-1\right){n \choose 2}\times n$, and $\left(n-1\right){n \choose 2}\times\left(n-1\right){n \choose 2}$.
Incidentally the matrices $\mathbf{U}$ and $\mathbf{V}$ denotes
$n{n \choose 2}\times n{n \choose 2}$ matrix which are respectively
associated with basis for the left and right eigenspaces respectively.
Finally the vectors $\boldsymbol{\mu}$, $\overline{\boldsymbol{\nu}}$
are $n{n \choose 2}$ dimensional vectors such that the entries of
their Hadamard product $\boldsymbol{\mu}\star\overline{\boldsymbol{\nu}}$
yield the eigenvalues of the adjacency matrix of $G_{2}$. The $k$-th
column vector of the matrices $\mathbf{U}\cdot\mbox{diag}\left\{ \boldsymbol{\mu}\right\} $
and $\overline{\mathbf{V}\cdot\mbox{diag}\left\{ \boldsymbol{\nu}\right\} }$
denote the scaled left eigenvectors and right scaled eigenvectors
respectively.\\
For $0\le j_{1}<j_{2}<n$ let the matrix $\mathbf{A}^{[j_{0},j_{1}]}$
denote $n\times n$ matrices defined by \textbf{
\[
\mathbf{A}^{\left[j_{1},j_{2}\right]}=\mathbf{A}\circ\sum_{\sigma\in S_{2}}\mathbf{e}_{\sigma\left(j_{1}\right)}\otimes\mathbf{e}_{\sigma\left(j_{2}\right)}^{T}
\]
}where the set $\left\{ \mathbf{e}_{j}\right\} _{0\le j<n}$ denotes
column vectors of the identity matrix $\mathbf{I}_{n}$. By construction
we have 
\[
\mathbf{A}=\sum_{0\le j_{1}<j_{2}<n}\mathbf{A}^{\left[j_{1},j_{2}\right]}.
\]
Furthermore let the spectral decomposition of the matrix $\mathbf{A}^{\left[j_{1},j_{2}\right]}$
be expressed as
\[
\mathbf{A}^{\left[j_{1},j_{2}\right]}=\left(\mathbf{U}^{\left[j_{1},j_{2}\right]}\cdot\mbox{diag}\left\{ \boldsymbol{\mu}^{\left[j_{1},j_{2}\right]}\right\} \right)\cdot\left(\mathbf{V}^{\left[j_{1},j_{2}\right]}\cdot\mbox{diag}\left\{ \boldsymbol{\nu}^{\left[j_{1},j_{2}\right]}\right\} \right)^{\dagger}
\]
\[
\mathbf{e}_{j_{1}}\cdot\mathbf{e}_{j_{1}}^{T}+\mathbf{e}_{j_{2}}\cdot\mathbf{e}_{j_{2}}^{T}=\mathbf{U}^{\left[j_{1},j_{2}\right]}\cdot\left(\mathbf{V}^{\left[j_{1},j_{2}\right]}\right)^{\dagger}.
\]
we have

\[
\forall\,0\le i_{1},i_{2}<n,\quad a_{i_{1},i_{2}}=
\]
\begin{equation}
\sum_{\begin{array}{c}
0\le k<n\\
0\le j_{1}<j_{2}<n
\end{array}}\left[\left(\sqrt{n-1}\,\mu_{k}^{\left[j_{1},j_{2}\right]}\right)\,\left(\frac{u_{i_{0},k}^{\left[j_{1},j_{2}\right]}}{\sqrt{n-1}}\right)\right]\overline{\left[\left(\sqrt{n-1}\,\nu_{k}^{\left[j_{1},j_{2}\right]}\right)\,\left(\frac{v_{i_{1},k}^{\left[j_{1},j_{2}\right]}}{\sqrt{n-1}}\right)\right]},\label{eq:Spectral_Decomposition}
\end{equation}
and
\[
\forall\,0\le i_{1},i_{2}<n,\quad\delta_{i_{1},i_{2}}=\sum_{\begin{array}{c}
0\le k<n\\
0\le j_{1}<j_{2}<n
\end{array}}\left(\frac{u_{i_{0},k}^{\left[j_{1},j_{2}\right]}}{\sqrt{n-1}}\right)\overline{\left(\frac{v_{i_{1},k}^{\left[j_{1},j_{2}\right]}}{\sqrt{n-1}}\right)}.
\]
The right-hand side of the expressions in \ref{eq:Spectral_Decomposition}
should be viewed as expressing inner-products of $n{n \choose 2}$-dimensional
vectors. The first $n$ rows of the matrices $\mathbf{U}$ and $\overline{\mathbf{V}}$
are therefore determined by the expansion above. The remaining $\left({n \choose 2}-1\right)n$
rows of the matrix $\mathbf{U}$ and $\overline{\mathbf{V}}$ are
determined by the Gram-Schmdit process.\\
The proposed construction used for adjacency matrices of graphs is
easily extended to higher order hypermatrices via the BM algebra.
For notational convenience we discuss here only the third order hypermartrice
case. We recall from \cite{GER} that by analogy to the matrix case
the spectral decomposition a third order hypermatrix $\mathbf{A}$
is expressed in terms of scaled eigenmatrices as follows 
\[
\mathbf{A}=\mbox{Prod}\left(\mbox{Prod}\left(\mathbf{Q},\mathbf{D}_{3},\mathbf{D}_{3}^{T}\right),\left[\mbox{Prod}\left(\mathbf{U},\mathbf{D}_{2},\mathbf{D}_{2}^{T}\right)\right]^{T^{2}},\left[\mbox{Prod}\left(\mathbf{V},\mathbf{D}_{1},\mathbf{D}_{1}^{T}\right)\right]^{T}\right).
\]
The collection of row-depth matrix slices of the hypermatrices $\mathbf{Q}$,
$\mathbf{U}$ and $\mathbf{V}$ yields bases for the \emph{eigenmatrices}
of $\mathbf{A}$ which are also subject to the non-correlation constraints
\[
\left[\mbox{Prod}\left(\mathbf{Q},\mathbf{U}^{T^{2}},\mathbf{V}^{T}\right)\right]_{i_{1},i_{2},i_{2}}=\begin{cases}
\begin{array}{cc}
1 & \mbox{ if }i_{1}=i_{2}=i_{3}\\
0 & \mbox{otherwise}
\end{array}\end{cases}.
\]
The scaling hypermatrices $\left\{ \mathbf{D}_{i}\right\} _{1\le i\le3}$
correspond to third order hypermatrix analog of diagonal matrices.
Consequently, the row-depth slices of the hypermatrices $\mbox{Prod}\left(\mathbf{Q},\mathbf{D}_{3},\mathbf{D}_{3}^{T}\right)$,
$\mbox{Prod}\left(\mathbf{U},\mathbf{D}_{2},\mathbf{D}_{2}^{T}\right)$,
and $\mbox{Prod}\left(\mathbf{V},\mathbf{D}_{1},\mathbf{D}_{1}^{T}\right)$,
correspond to scaled eigenmatrices of $\mathbf{A}$. By analogy to
the matrix case, for all triplets label $\left(j_{1},j_{2},j_{3}\right)$
for which the inequality $0\le j_{1}<j_{2}<j_{3}<n$ is satisfied,
we define the matrix $\mathbf{A}^{[j_{1},j_{2},j_{3}]}$ to be the
third order hypermatrix expressed by \textbf{
\begin{equation}
\mathbf{A}^{[j_{1},j_{2},j_{3}]}=\mathbf{A}\circ\sum_{\sigma\in S_{3}}\mathbf{e}_{\sigma\left(j_{1}\right)}\otimes\mathbf{e}_{\sigma\left(j_{2}\right)}^{T}\otimes\mathbf{e}_{\sigma\left(j_{3}\right)}^{T^{2}}.
\end{equation}
}where the set $\left\{ \mathbf{e}_{j}\right\} _{0\le j<n}$ denotes
column vectors of the identity matrix $\mathbf{I}_{n}$ and the transpose
operation here refers to the cyclic permutation of the hypermatric
entries as introduced in \cite{GER}. Similarly we have 
\begin{equation}
\mathbf{A}=\sum_{0\le j_{1}<j_{2}<j_{3}<n}\mathbf{A}^{\left[j_{1},j_{2},j_{3}\right]}
\end{equation}
by appropriately concatenating the spectral decomposition of the hypermatrix
minors $\mathbf{A}^{[j_{0},j_{1},j_{2}]}$ and using the generalization
to hypermatrices of the constrained inverse matrix pair problem, we
deduce the spectral decomposition of a larger $n{n \choose 3}\times n{n \choose 3}\times n{n \choose 3}$
which admits the hypermatrix $\mathbf{A}$ as a sub-hypermatrix.\\
More generally a similar construction is easily devised for higher
order hypermatrices. The argument therefore provides a way to reduce
the spectral decomposition of $k$-th order cubic hypermatrices to
the spectral decomposition of cubic sub-hypermatrices with sides of
length $k$.
\end{proof}

\section{Computing spectral elimination ideals.}

For the purpose of introducing new iterative procedures to approximate
the spectral decomposition of matrices and hypermartices, We describe
here two elimination procedures for computing generators of elimination
ideals associated with matrix and hypermatrix spectral constraints.
The first elimination method uses only properties of the algebra of
matrices to obtain generators for the ideal associated with eigenvalues
in the matrix case and scaling values in the general case of higher
order hypermatrices. The first elimination method also allows us to
derive new generalizations of the matrix determinant polynomial for
higher order hypermatrices. The second elimination method uses a Hypermatrix
generalization of the classical Paserval identity to obtain generators
for the ideal associated with entries of uncorrelated tuples (which
are hypermatrix analog of matrix eigenvectors). The iterative procedure
described subsequently will use the generators devised by the second
method.

\subsection{Spectral elimination ideals for matrices.}

We start by discussing both elimination method in the matrix case
and subsequently proceed to extend the arguments to hypermatrices
of arbitrary orders. Let us recall for convenience of the reader the
well known matrix spectral decomposition of an $n\times n$ matrix
$\mathbf{A}$. Such a decomposition is obtained by solving for $n\times n$
matrices $\mathbf{U}$, $\mathbf{V}$, and diagonal matrices $\left\{ \mathbf{D}_{i}\right\} _{1\le i\le2}$
subject to the spectral constraints 
\begin{equation}
\begin{cases}
\begin{array}{ccc}
\mathbf{A} & = & \left(\mathbf{U}\cdot\mathbf{D}_{1}\right)\cdot\left(\overline{\mathbf{V}\cdot\mathbf{D}_{2}}\right)^{T}\\
\mathbf{U}\cdot\overline{\mathbf{V}}^{T} & = & \mathbf{I}_{n}\\
\mathbf{D}_{i}\circ\mathbf{D}_{i} & = & \mathbf{D}_{i}^{T}\cdot\mathbf{D}_{i},\quad\forall\:0\le i<2
\end{array}\end{cases}.\label{eq:Spectral_Constraints}
\end{equation}
For notational convenience we reformulate the diagonality constraints
in the spectral decomposition, in terms of $n$-dimensional vectors
$\boldsymbol{\mu}$ and $\boldsymbol{\nu}$ as follows 
\[
\begin{cases}
\begin{array}{ccc}
\mathbf{A} & = & \left(\mathbf{U}\cdot\mbox{diag}\left\{ \boldsymbol{\mu}\right\} \right)\cdot\left(\mathbf{V}\cdot\mbox{diag}\left\{ \boldsymbol{\nu}\right\} \right)^{\dagger}\\
\mathbf{U}\cdot\mathbf{V}^{\dagger} & = & \mathbf{I}_{n}
\end{array}\end{cases}.
\]
The two elimination ideals obtained by the two distinct method derived
from the spectral decomposition constraints are : 
\[
\mathcal{I}_{\boldsymbol{\mu}\circ\overline{\boldsymbol{\nu}}}\,:=\mathbb{C}\left[\boldsymbol{\mu}\circ\overline{\boldsymbol{\nu}}\right]\cap\mbox{ Ideal generated by}\left\{ \left(\mathbf{U}\cdot\mbox{diag}\left\{ \boldsymbol{\mu}\right\} \right)\cdot\left(\mathbf{V}\cdot\mbox{diag}\left\{ \boldsymbol{\nu}\right\} \right)^{\dagger}-\mathbf{A},\:\mathbf{U}\cdot\mathbf{V}^{\dagger}-\mathbf{I}_{n}\right\} ,
\]
\[
\mbox{and}
\]
\[
\mathcal{I}_{\mathbf{U},\overline{\mathbf{V}}}\,:=\mathbb{C}\left[\mathbf{U},\overline{\mathbf{V}}\right]\cap\mbox{ Ideal generated by}\left\{ \left(\mathbf{U}\cdot\mbox{diag}\left\{ \boldsymbol{\mu}\right\} \right)\cdot\left(\mathbf{V}\cdot\mbox{diag}\left\{ \boldsymbol{\nu}\right\} \right)^{\dagger}-\mathbf{A},\,\mathbf{U}\cdot\mathbf{V}^{\dagger}-\mathbf{I}_{n}\right\} .
\]
Let us start by describing the derivation of generators for $\mathcal{I}_{\boldsymbol{\mu}\circ\overline{\boldsymbol{\nu}}}$.
The starting point for the derivation is the matrix identity 
\[
\forall\:0\le k\le n,\quad\mathbf{A}^{k}=\left(\mathbf{U}\cdot\mbox{diag}\left\{ \boldsymbol{\mu}^{\circ^{k}}\right\} \right)\cdot\left(\mathbf{V}\cdot\mbox{diag}\left\{ \boldsymbol{\nu}^{\circ^{k}}\right\} \right)^{\dagger}.
\]
Let $\left\{ \mathbf{u}_{i}\right\} _{0\le i<n}$ and $\left\{ \overline{\mathbf{v}_{j}}\right\} _{0\le j<n}$
denote respectively row vectors of the matrices $\mathbf{U}$ and
$\overline{\mathbf{V}}$, we reformulate these identities into the
following Vandermonde type equalities of the form 
\[
\forall\,0\le i,j<n,\quad\mathbf{u}_{i}\circ\overline{\mathbf{v}_{j}}=\left(\mbox{Vandermonde}\left\{ \boldsymbol{\mu}\circ\overline{\boldsymbol{\nu}}\right\} \right)^{-1}\cdot\left(\begin{array}{c}
\left[\mathbf{A}^{0}\right]_{i,j}\\
\vdots\\
\left[\mathbf{A}^{n-1}\right]_{i,j}
\end{array}\right)
\]
the constraints can be combined to yield the following identity. 
\[
\left(\begin{array}{c}
\left[\mathbf{u}_{i}\circ\overline{\mathbf{v}_{j}}\right]_{0}\\
\vdots\\
\left[\mathbf{u}_{i}\circ\overline{\mathbf{v}_{j}}\right]_{n-1}
\end{array}\right)_{0\le i,j<n}=\left(\mathbf{I}_{n^{2}}\otimes\mbox{Vandermonde}\left(\boldsymbol{\mu}\circ\overline{\boldsymbol{\nu}}\right)\right)^{-1}\cdot\left(\begin{array}{c}
\left[\mathbf{A}^{0}\right]_{i,j}\\
\vdots\\
\left[\mathbf{A}^{n-1}\right]_{i,j}
\end{array}\right)_{0\le i,j<n}.
\]
It is implicitly assumed in the identity above that $\mathbf{A}$
has distinct eigenvalues, otherwise the identity above is expressed
using the Penrose matrix inverse instead as follows 
\[
\left(\begin{array}{c}
\left[\mathbf{u}_{i}\circ\overline{\mathbf{v}_{j}}\right]_{0}\\
\vdots\\
\left[\mathbf{u}_{i}\circ\overline{\mathbf{v}_{j}}\right]_{n-1}
\end{array}\right)_{0\le i,j<n}=\left(\mathbf{I}_{n^{2}}\otimes\mbox{Vandermonde}\left(\boldsymbol{\mu}\circ\overline{\boldsymbol{\nu}}\right)\right)^{+}\cdot\left(\begin{array}{c}
\left[\mathbf{A}^{0}\right]_{i,j}\\
\vdots\\
\left[\mathbf{A}^{n-1}\right]_{i,j}
\end{array}\right)_{0\le i,j<n}.
\]
Having thus expressed in the equality above the entries of the vectors
$\left\{ \mathbf{u}_{i}\circ\overline{\mathbf{v}_{j}}\right\} _{0\le i,j<n}$
only in terms of the entries of $\mathbf{A}$ and its eigenvalues,
we derive the generators for $\mathcal{I}_{\boldsymbol{\mu}\circ\overline{\boldsymbol{\nu}}}$
as prescribed by tautologies 
\[
\mathcal{I}_{\boldsymbol{\mu}\circ\overline{\boldsymbol{\nu}}}=\mbox{Ideal generated by}\left\{ \left(\mathbf{u}_{i}\circ\overline{\mathbf{v}_{j}}\right)\circ\left(\mathbf{u}_{j}\circ\overline{\mathbf{v}_{i}}\right)-\left(\mathbf{u}_{i}\circ\overline{\mathbf{v}_{i}}\right)\circ\left(\mathbf{u}_{j}\circ\overline{\mathbf{v}_{j}}\right)\right\} _{0\le i,j<n}.
\]
Note that the computation of resultants or alternatively the computation
of Groebner basis can also be used to compute generators for $\mathcal{I}_{\boldsymbol{\mu}\circ\overline{\boldsymbol{\nu}}}$
starting from the spectral constraints in \ref{eq:Spectral_Constraints}
as pointed out in \cite{GER}. However these approaches are less direct
and considerably less efficient because one must identify for spectral
constraints good monomial orderings. The elimination method also has
the benefit of producing explicit expressions for the characteristic
polynomial and determinant. We illustrate this by considering the
case $n=2$, where there generator of the ideal $\mathcal{I}_{\boldsymbol{\mu}\circ\overline{\boldsymbol{\nu}}}$
results in a single vector equality 
\[
\left(\mathbf{u}_{0}\circ\overline{\mathbf{v}_{1}}\right)\circ\left(\mathbf{u}_{1}\circ\overline{\mathbf{v}_{0}}\right)-\left(\mathbf{u}_{0}\circ\overline{\mathbf{v}_{0}}\right)\circ\left(\mathbf{u}_{1}\circ\overline{\mathbf{v}_{1}}\right)=\mathbf{0}_{2\times1}
\]
and yields the equality
\[
\left(\begin{array}{c}
\left(\mu_{1}\overline{\nu_{1}}\right)^{2}-\left(a_{0,0}+a_{1,1}\right)\mu_{1}\overline{\nu_{1}}+\left(a_{0,0}a_{1,1}-a_{0,1}a_{1,0}\right)\\
\left(\mu_{0}\overline{\nu_{0}}\right)^{2}-\left(a_{0,0}+a_{1,1}\right)\mu_{0}\overline{\nu_{0}}+\left(a_{0,0}a_{1,1}-a_{0,1}a_{1,0}\right)
\end{array}\right)=\left(\begin{array}{c}
0\\
0
\end{array}\right).
\]

Let us now turn our attention to the derivation of generators for
the elimination ideal $\mathcal{I}_{\mathbf{U},\overline{\mathbf{V}}}$.
The derivation is achieved by eliminating the variables associated
with the eigenvalues determined by the entries of the vectors $\boldsymbol{\mu}$
and $\overline{\boldsymbol{\nu}}$ via Parseval's identity. For the
purposes of the proposed elimination method we view the spectral constraints
in \ref{eq:Spectral_Constraints} as a collection of $n^{2}$ inner
product equalities of the form 
\[
\left\{ \left\langle \mathbf{u}_{i_{0}}\circ\boldsymbol{\mu},\,\overline{\boldsymbol{\nu}\circ\mathbf{v}_{i_{1}}}\right\rangle =a_{i_{0},i_{1}}\right\} _{0\le i_{0},i_{1}<n},
\]
where the sets $\left\{ \mathbf{u}_{i}\right\} _{0\le i<n}$ and $\left\{ \mathbf{v}_{i}\right\} _{0\le i<n}$
denote row vectors of the matrix $\mathbf{U}$ and $\mathbf{V}$ respectively.
By Parsevals' identity the $n^{2}$ constraints can be reformulated
as 
\begin{equation}
\left\{ \sum_{0\le k<n}\left\langle \mathbf{u}_{i_{0}}\circ\boldsymbol{\mu},\overline{\boldsymbol{\nu}\circ\mathbf{v}_{i_{1}}}\right\rangle _{\left(\overline{\mathbf{v}}_{k}\right)^{T}\cdot\mathbf{u}_{k}}=a_{i,j}=\sum_{0\le j_{0},j_{1}<n}\mu_{j_{0}}\overline{\nu_{j_{1}}}\:f_{n\cdot j_{0}+j_{1},n\cdot i_{0}+i_{1}}\left(\mathbf{U},\overline{\mathbf{V}}\right)\right\} _{0\le i_{0},i_{1}<n}\label{eq:Resolution_of_identity}
\end{equation}
where $\left\{ f_{n\cdot j_{0}+j_{1},n\cdot i_{0}+i_{1}}\left(\mathbf{U},\overline{\mathbf{V}}\right)\right\} _{0\le n\cdot j_{0}+j_{1},n\cdot i_{0}+i_{1}<n^{2}}\subset\mathbb{C}\left[\mathbf{U},\overline{\mathbf{V}}\right]$
and expressed by 
\[
f_{n\cdot j_{0}+j_{1},n\cdot i_{0}+i_{1}}\left(\mathbf{U},\overline{\mathbf{V}}\right)=u_{i_{0}j_{0}}\,\left(\sum_{0\le k<n}\overline{v_{j_{0}k}}u_{j_{1}k}\right)\,\overline{v_{i_{1}j_{1}}}.
\]
The constraints form a system of $n^{2}$ equations in the $n^{2}$
unknowns $\left\{ \mu_{i}\,\overline{\nu_{j}}\right\} _{0\le i,j<n}$.
We therefore express the polynomial constraints in \ref{eq:Resolution_of_identity}
as follows 
\begin{equation}
\left(\begin{array}{ccccc}
f_{0,0} & \cdots & f_{0,j} & \cdots & f_{0,\left(n^{2}-1\right)}\\
\vdots & \ddots & \vdots & \iddots & \vdots\\
f_{i,0} & \cdots & f_{ij} & \cdots & f_{i,\left(n^{2}-1\right)}\\
\vdots & \iddots & \vdots & \ddots & \vdots\\
f_{\left(n^{2}-1\right),0} & \cdots & f_{\left(n^{2}-1\right),j} & \cdots & f_{\left(n^{2}-1\right),\left(n^{2}-1\right)}
\end{array}\right)\left(\begin{array}{c}
\mu_{0}\overline{\nu_{0}}\\
\vdots\\
\mu_{i}\overline{\nu_{j}}\\
\vdots\\
\mu_{n-1}\overline{\nu_{n-1}}
\end{array}\right)=\left(\begin{array}{c}
a_{0,0}\\
\vdots\\
a_{i,j}\\
\vdots\\
a_{\left(n-1\right),\left(n-1\right)}
\end{array}\right).
\end{equation}
For any integer $0\le k<n^{2}$, let $\mathbf{F}_{k}$ denote the
$n^{2}\times n^{2}$ matrix constructed by substituting the $k$-th
column of the left hand side matrix $\mathbf{F}$ above by the righthand
side $n^{2}\times1$ vector made up by the entries of $\mathbf{A}$
as prescribed by the classical Cramer's rule. The generators for $I_{\mathbf{U},\overline{\mathbf{V}}}$
are determined by rational function identities derived from tautologies
similar to the tautologies used in the first elimination method. 
\[
\left\{ \left(\mu_{i}\overline{\nu_{i}}\right)\left(\mu_{j}\overline{\nu_{j}}\right)=\left(\mu_{i}\overline{\nu_{j}}\right)\left(\mu_{j}\overline{\nu_{i}}\right)\right\} _{0\le i<j<n}.
\]
It follows from Cramer's rule that generators for the ideal $\mathcal{I}_{\mathbf{U},\overline{\mathbf{V}}}$
are derive from the rational identities 
\[
\left\{ \frac{\det\left(\mathbf{F}_{n\cdot i+i}\cdot\mathbf{F}_{n\cdot j+j}\right)-\det\left(\mathbf{F}_{n\cdot i+j}\cdot\mathbf{F}_{n\cdot j+i}\right)}{\det\left(\mathbf{F}^{2}\right)}=0\right\} _{0\le i<j<n}.
\]
\medskip{}
The elimination ideal $\mathcal{I}_{\mathbf{U},\overline{\mathbf{V}}}$
is seldom used in the literature however it's analog is of crucial
importance for approximating the spectrum of hypermatrices.

\subsection{Spectral elimination ideals for hypermatrices.}

For notational convenience we restrict the discussion here to third
order hypermatrices. We recall the spectral decomposition for third
order hypermatrices is expressed by 
\[
\mathbf{A}=\mbox{Prod}\left(\mbox{Prod}\left(\mathbf{Q},\mathbf{D}_{3},\mathbf{D}_{3}^{T}\right),\left[\mbox{Prod}\left(\mathbf{U},\mathbf{D}_{2},\mathbf{D}_{2}^{T}\right)\right]^{T^{2}},\left[\mbox{Prod}\left(\mathbf{V},\mathbf{D}_{1},\mathbf{D}_{1}^{T}\right)\right]^{T}\right).
\]
The collection of matrix slices of the hypermatrices $\mathbf{Q}$,
$\mathbf{U}$ and $\mathbf{V}$ collects the \emph{eigen-matrices}
of $\mathbf{A}$ which are subject to the constraints
\[
\left[\mbox{Prod}\left(\mathbf{Q},\,\mathbf{U}^{T^{2}},\,\mathbf{V}^{T}\right)\right]_{i_{1},i_{2},i_{2}}=\boldsymbol{\Delta}=\begin{cases}
\begin{array}{cc}
1 & \mbox{ if }i_{1}=i_{2}=i_{3}\\
0 & \mbox{otherwise}
\end{array}\end{cases}.
\]
The scaling hypermatrices $\left\{ \mathbf{D}_{i}\right\} _{1\le i\le3}$
are hypermatrix analog of diagonal matrices. Third order hypermatrix
diagonality constraints are similar to matrix diagonality constraints
and expressed by 
\[
\forall\,1\le i\le3,\quad\mathbf{D}_{i}\circ\mathbf{D}_{i}\circ\mathbf{D}_{i}=\mbox{Prod}\left(\mathbf{D}_{i}^{T},\,\mathbf{D}_{i}^{T^{2}},\,\mathbf{D}_{i}\right),
\]
The slices of $\mbox{Prod}\left(\mathbf{Q},\mathbf{D}_{3},\mathbf{D}_{3}^{T}\right)$,
$\mbox{Prod}\left(\mathbf{U},\mathbf{D}_{2},\mathbf{D}_{2}^{T}\right)$,
and $\mbox{Prod}\left(\mathbf{V},\mathbf{D}_{1},\mathbf{D}_{1}^{T}\right)$,
correspond to scaled eigenmatrices of $\mathbf{A}$. Moreover the
spectral decomposition constraints can rewritten in terms of inner-product
constraints of the form 
\[
\forall\:0\le i,j,k<2,\quad a_{ijk}=\left\langle \left(\boldsymbol{\alpha}_{i}\circ\mathbf{q}_{ik}\circ\boldsymbol{\alpha}_{k}\right),\,\left(\boldsymbol{\beta}_{j}\circ\mathbf{u}_{ji}\circ\boldsymbol{\beta}_{i}\right),\,\left(\boldsymbol{\gamma}_{k}\circ\mathbf{v}_{kj}\circ\boldsymbol{\gamma}_{j}\right)\right\rangle 
\]
\[
\left\langle \mathbf{q}_{ik},\,\mathbf{u}_{ji},\,\mathbf{v}_{kj}\right\rangle =\delta_{i,j,k}=\begin{cases}
\begin{array}{cc}
1 & \mbox{if }i=j=k\\
0 & \mbox{otherwise }
\end{array} & \forall\:0\le i,j,k<n\end{cases}.
\]
The derivation according to the first elimination method is best illustrated
for cubic hypermatrices having side length equal to $2$. In particular
for $2\times2\times2$ the spectral decomposition yields the following
Vandermonde type equalities generally expressed 
\[
\forall\:0\le i,j,k<2,\quad\mathbf{q}_{ik}\circ\mathbf{u}_{ji}\circ\mathbf{v}_{kj}=
\]
\[
\left(\mbox{Vandermonde}\left\{ \left(\boldsymbol{\alpha}_{i}\circ\boldsymbol{\alpha}_{k}\right)\circ\left(\boldsymbol{\beta}_{j}\circ\boldsymbol{\beta}_{i}\right)\circ\left(\boldsymbol{\gamma}_{k}\circ\boldsymbol{\gamma}_{j}\right)\right\} \right)^{-1}\left(\begin{array}{c}
\delta_{i,j,k}\\
a_{i,j,k}
\end{array}\right)
\]
Having thus expressed in the equalities above the entries of the vectors
$\left\{ \mathbf{q}_{ik}\circ\mathbf{u}_{ji}\circ\mathbf{v}_{kj}\right\} _{0\le i,j,k<n}$
only in terms of the entries of $\mathbf{A}$ and its scaling values,
we derive the generators for $\mathcal{I}_{\boldsymbol{\alpha},\boldsymbol{\beta},\boldsymbol{\gamma}}$
as prescribed by the tautology
\[
\left(\mathbf{q}_{00}\circ\mathbf{u}_{00}\circ\mathbf{v}_{00}\right)\circ\left(\mathbf{q}_{01}\circ\mathbf{u}_{10}\circ\mathbf{v}_{11}\right)\circ\left(\mathbf{q}_{11}\circ\mathbf{u}_{01}\circ\mathbf{v}_{10}\right)\circ\left(\mathbf{q}_{10}\circ\mathbf{u}_{11}\circ\mathbf{v}_{01}\right)-
\]
\[
\left(\mathbf{q}_{01}\circ\mathbf{u}_{00}\circ\mathbf{v}_{10}\right)\circ\left(\mathbf{q}_{00}\circ\mathbf{u}_{10}\circ\mathbf{v}_{01}\right)\circ\left(\mathbf{q}_{10}\circ\mathbf{u}_{01}\circ\mathbf{v}_{00}\right)\circ\left(\mathbf{q}_{11}\circ\mathbf{u}_{11}\circ\mathbf{v}_{11}\right)=\mathbf{0}_{2\times1}
\]
which results in third order hypermatrix analog of $2\times2\times2$
characteristic polynomial 
\[
\left(\begin{array}{c}
a_{001}a_{010}a_{100}\alpha_{11}^{2}\beta_{11}^{2}\gamma_{11}^{2}-a_{011}a_{101}a_{110}\alpha_{01}^{2}\beta_{01}^{2}\gamma_{01}^{2}+a_{000}a_{011}a_{101}a_{110}-a_{001}a_{010}a_{100}a_{111}\\
a_{001}a_{010}a_{100}\alpha_{01}^{2}\beta_{01}^{2}\gamma_{01}^{2}-a_{011}a_{101}a_{110}\alpha_{00}^{2}\beta_{00}^{2}\gamma_{00}^{2}+a_{000}a_{011}a_{101}a_{110}-a_{001}a_{010}a_{100}a_{111}
\end{array}\right)
\]
\[
=\left(\begin{array}{c}
0\\
0
\end{array}\right).
\]
The first elimination method therefore yield hypermatrix analog of
characteristic polynomials as well as an expression of the hyperdeterminants
for hypermatrices of size $2^{m}\times2^{m}\times2^{m}$ which arises
as Kronecker products of $2\times2\times2$ hypermatrices. The first
method is considerably more difficult to extend to arbitrary hypermatrices
of size $n\times n\times n$ because of the lack of Vandermonde for
$n>2$. However, the second elimination method extend to hypermatrices
of size $n\times n\times n$. Similarly to the matrix case, the second
elimination method based on a hypermatrix generalization of Parseval's
identity. We eliminate the variables associated with the scaling values
by considering the sequence of hypermatrix defined the recurrence
relations 
\[
\mathbf{G}_{0}=\boldsymbol{\Delta}_{2},\quad\mathbf{G}_{k+1}=\mbox{Prod}_{\mathbf{G}_{k}}\left(\mathbf{Q},\mathbf{U}^{T^{2}},\mathbf{V}^{T}\right),
\]
in conjunction with the constraints 
\[
\left\{ \mathbf{A}=\mbox{Prod}_{\mathbf{G}_{k}}\left(\mbox{Prod}\left(\mathbf{Q},\mathbf{D}_{0},\mathbf{D}_{0}^{T}\right),\,\left[\mbox{Prod}\left(\mathbf{U},\mathbf{D}_{2},\mathbf{D}_{2}^{T}\right)\right]^{T^{2}},\,\left[\mbox{Prod}\left(\mathbf{V},\mathbf{D}_{1},\mathbf{D}_{1}^{T}\right)\right]^{T}\right)\right\} _{k<n}
\]
Finally using Cramer's rule we express the the scaling variables as
rational functions of the entries of the $\mathbf{Q}$,$\mathbf{U}$,$\mathbf{V}$
and $\mathbf{A}$. As illustration for the second elimination method
we consider the case of a $2\times2\times2$ hypermatrix $\mathbf{A}$
such that $\mathbf{A}^{T}=\mathbf{A}$, which admit a spectral decomposition
where 
\[
\mathbf{D}_{0}=\mathbf{D}_{1}=\mathbf{D}_{2}\mbox{ and }\mathbf{Q}=\mathbf{U}=\mathbf{V}.
\]
Using Cramer's rule to isolate the scaling variables from the constraints
\[
\left\{ \mathbf{A}=\mbox{Prod}_{\mathbf{G}_{k}}\left(\mbox{Prod}\left(\mathbf{Q},\mathbf{D}_{0},\mathbf{D}_{0}^{T}\right),\left[\mbox{Prod}\left(\mathbf{U},\mathbf{D}_{2},\mathbf{D}_{2}^{T}\right)\right]^{T^{2}},\left[\mbox{Prod}\left(\mathbf{V},\mathbf{D}_{1},\mathbf{D}_{1}^{T}\right)\right]^{T}\right)\right\} _{k<2}
\]
we derive the elimination ideal $\mathcal{I}_{\mathbf{Q}}$ from the
rational identities prescribed by the tautologies :
\[
\begin{cases}
\begin{array}{ccc}
\left(\lambda_{km}^{4}\lambda_{kp}^{2}\right)^{3} & = & \left[\left(\lambda_{km}^{2}\right)^{3}\right]^{2}\left[\left(\lambda_{kp}^{2}\right)^{3}\right]\\
\left(\lambda_{km}^{2}\lambda_{kn}^{2}\lambda_{kp}^{2}\right)^{3} & = & \left(\lambda_{km}^{2}\right)^{3}\left(\lambda_{kn}^{2}\right)^{3}\left(\lambda_{kp}^{2}\right)^{3}
\end{array} & .\end{cases}
\]
Unlike the first elimination method, the second elimination method
has no restriction on the size of the hypermatrices. Furthermore,
for the general purpose of approximating the spectrum of arbitrary
hypermatrices subsequent methods discussed here construct approximation
arbitrary size hypermatrices from spectral decomposition of the $2\times2\times2$
minors.

\section{\label{sec:Approximating}Approximating the spectral decompositions
of matrices and hypermatrices.}

We describe here a recursive construction for approximating the spectral
decomposition of matrices and hypermatrices based on a refinement
of the proof of theorem \ref{thm:Spectral decomposition}. We start
by discussing the matrix case and subsequently briefly discuss how
the arguments are extended to hypermatrices. For some $n\times n$
matrix $\mathbf{A}$ with complex entries, let $\mathbf{A}_{\tau}$
denote the matrix minor constructed as follows
\[
\mathbf{A}_{\tau}=
\]
\[
\mathbf{A}\circ\left[\frac{1}{n-2}\left(\left(\sum_{0\le i\ne\tau<n}\mathbf{e}_{i}\right)\cdot\left(\sum_{0\le j\ne\tau<n}\mathbf{e}_{j}\right)^{T}-\sum_{0\le k\ne\tau<n}\mathbf{e}_{k}\cdot\mathbf{e}_{k}^{T}\right)+\frac{1}{n-1}\sum_{0\le k\ne\tau<n}\mathbf{e}_{k}\cdot\mathbf{e}_{k}^{T}\right].
\]
The minors are constructed as to obtain the identity 
\[
\mathbf{A}=\sum_{0\le\tau<n}\mathbf{A}_{\tau}.
\]
We further assume for the sake of the argument that the spectral decomposition
of the matrix minors $\left\{ \mathbf{A}_{\tau}\right\} _{0\le\tau<n}$
are known and given by 
\[
\mathbf{A}_{\tau}=\left(\mathbf{U}^{[\tau]}\cdot\mbox{diag}\left\{ \boldsymbol{\mu}^{[\tau]}\right\} \right)\cdot\left(\mathbf{V}^{[\tau]}\cdot\mbox{diag}\left\{ \boldsymbol{\nu}^{[\tau]}\right\} \right)^{\dagger},\quad\mathbf{I}_{n}-\mathbf{e}_{\tau}\cdot\mathbf{e}_{\tau}^{T}=\mathbf{U}^{[\tau]}\cdot\left(\mathbf{V}^{[\tau]}\right)^{\dagger}.
\]
By concatenating the spectral decomposition of the matrices $\left\{ \mathbf{A}_{\tau}\right\} _{0\le\tau<n}$
as was done in the proof of \ref{thm:Spectral decomposition}, we
obtain that for all $0\le i_{0},i_{1}<n$ 
\[
a_{i_{0},i_{1}}=\sum_{0\le\tau,t<n}\left[\left(\sqrt{n-1}\,\mu_{t}^{\left[\tau\right]}\right)\,\left(\frac{u_{i_{0},t}^{\left[\tau\right]}}{\sqrt{n-1}}\right)\right]\overline{\left[\left(\sqrt{n-1}\,\nu_{t}^{\left[\tau\right]}\right)\,\left(\frac{v_{i_{1},t}^{\left[\tau\right]}}{\sqrt{n-1}}\right)\right]}
\]
\[
\delta_{i_{0},i_{1}}=\sum_{0\le\tau,t<n}\left(\frac{u_{i_{0}t}^{\left[\tau\right]}}{\sqrt{n-1}}\right)\overline{\left(\frac{v_{i_{1}t}^{\left[\tau\right]}}{\sqrt{n-1}}\right)}
\]
As was done in the proof of theorem \ref{thm:Spectral decomposition}
we construct the spectral decomposition of a larger $n^{2}\times n^{2}$
matrix (by concatenating the spectral decomposition of the matrices
$\left\{ \mathbf{A}_{\tau}\right\} _{0\le\tau<n}$) which admits $\mathbf{A}$
as a sub-matrix and expressed as
\[
\left(\begin{array}{cc}
\mathbf{A} & \mathbf{B}_{01}\\
\mathbf{B}_{10} & \mathbf{B}_{11}
\end{array}\right)=\left(\mathbf{U}\cdot\mbox{diag}\left\{ \boldsymbol{\mu}\right\} \right)\cdot\left(\mathbf{V}\cdot\mbox{diag}\left\{ \boldsymbol{\nu}\right\} \right)^{\dagger},\;\mathbf{I}_{n^{2}}=\mathbf{U}\cdot\mathbf{V}^{\dagger}.
\]
Following this construction, we discuss two approximation techniques
which enable a recursive approximation for the spectrum of the original
matrix. The first method is the inflation approximation technique.
It iteratively modifies the spectrum of the larger $n^{2}\times n^{2}$
matrices, in order to approximate the spectrum of the smaller $n\times n$
matrix. We start from the matrices $\mathbf{U}$ and $\mathbf{V}$
iteratively attempts to converge via gradient descent to matrices
\[
\left(\begin{array}{cc}
\mathbf{U}^{\prime} & \mathbf{0}_{n\times n\left(n-1\right)}\\
\mathbf{0}_{n\left(n-1\right)\times n} & \mathbf{Q}
\end{array}\right),\mbox{ and }\:\left(\begin{array}{cc}
\mathbf{V}^{\prime} & \mathbf{0}_{n\times n\left(n-1\right)}\\
\mathbf{0}_{n\left(n-1\right)\times n} & \left(\mathbf{Q}^{-1}\right)^{\dagger}
\end{array}\right)
\]
for which we have 
\[
\mathbf{I}_{n^{2}}=\left(\begin{array}{cc}
\mathbf{U}^{\prime} & \mathbf{0}_{n\times n\left(n-1\right)}\\
\mathbf{0}_{n\left(n-1\right)\times n} & \mathbf{Q}
\end{array}\right)\cdot\left(\begin{array}{cc}
\left(\mathbf{V}^{\prime}\right)^{\dagger} & \mathbf{0}_{n\times n\left(n-1\right)}\\
\mathbf{0}_{n\left(n-1\right)\times n} & \mathbf{Q}^{-1}
\end{array}\right)
\]
and
\[
\left(\begin{array}{cc}
\mathbf{A} & \mathbf{0}_{n\times n\left(n-1\right)}\\
\mathbf{0}_{n\left(n-1\right)\times n} & \mathbf{I}_{n\left(n-1\right)}
\end{array}\right)=
\]
\[
\left(\begin{array}{cc}
\mathbf{U}^{\prime}\cdot\mbox{diag}\left\{ \boldsymbol{\mu}^{\prime}\right\}  & \mathbf{0}_{n\times n\left(n-1\right)}\\
\mathbf{0}_{n\left(n-1\right)\times n} & \mathbf{Q}
\end{array}\right).\left(\begin{array}{cc}
\left(\mathbf{V}^{\prime}\cdot\mbox{diag}\left\{ \boldsymbol{\nu}^{\prime}\right\} \right)^{\dagger} & \mathbf{0}_{n\times n\left(n-1\right)}\\
\mathbf{0}_{n\left(n-1\right)\times n} & \mathbf{Q}^{-1}
\end{array}\right).
\]
We therefore derive from the resulting expansion a spectral decomposition
for $\mathbf{A}$ expressed by 
\begin{equation}
\mathbf{A}=\left(\mathbf{U}^{\prime}\cdot\mbox{diag}\left\{ \boldsymbol{\mu}^{\prime}\right\} \right)\cdot\left(\mathbf{V}^{\prime}\cdot\mbox{diag}\left\{ \boldsymbol{\nu}^{\prime}\right\} \right)^{\dagger}
\end{equation}
In contrast for large enough sizes we devise instead a second approximation
technique referred to as the truncation method. This method starts
from the spectral decomposition of the larger $n^{2}\times n^{2}$
matrix expressed by 
\begin{equation}
\left(\begin{array}{cc}
\mathbf{A} & \mathbf{B}_{01}\\
\mathbf{B}_{10} & \mathbf{B}_{11}
\end{array}\right)=\left(\mathbf{U}\cdot\mbox{diag}\left\{ \boldsymbol{\mu}\right\} \right)\cdot\left(\mathbf{V}\cdot\mbox{diag}\left\{ \boldsymbol{\nu}\right\} \right)^{\dagger},\;\mathbf{I}_{n^{2}}=\mathbf{U}\cdot\mathbf{V}^{\dagger}.
\end{equation}
where the vectors $\left\{ \mathbf{u}_{i}\right\} _{0\le i<n^{2}}$
and $\left\{ \overline{\mathbf{v}_{j}}\right\} _{0\le j<n^{2}}$ denote
$n^{2}$-dimensional column vectors of the matrices $\mathbf{U}$
and $\overline{\mathbf{V}}$ respectively. The approximation is therefore
obtained by truncating the vectors in the matrix to obtain 
\[
\left(\begin{array}{cc}
\widetilde{\mathbf{A}} & \mathbf{0}_{n\times\left(n^{2}-n\right)}\\
\mathbf{0}_{\left(n^{2}-n\right)\times n} & \mathbf{0}_{\left(n^{2}-n\right)\times\left(n^{2}-n\right)}
\end{array}\right)=
\]
\[
\sum_{0\le k<n}\left(\mu_{k}\mathbf{u}_{k}\circ\left(\begin{array}{c}
\mathbf{1}_{n\times1}\\
\mathbf{0}_{\left(n^{2}-n\right)\times1}
\end{array}\right)\right)\cdot\left(\nu_{k}\mathbf{v}_{k}\circ\left(\begin{array}{c}
\mathbf{1}_{n\times1}\\
\mathbf{0}_{\left(n^{2}-n\right)\times1}
\end{array}\right)\right)^{\dagger}.
\]
The total error incurred by the truncation from the original larger
matrix is
\[
\left\Vert \left(\begin{array}{cc}
\mathbf{A} & \mathbf{B}_{01}\\
\mathbf{B}_{01}^{\dagger} & \mathbf{B}_{11}
\end{array}\right)-\left(\begin{array}{cc}
\widetilde{\mathbf{A}} & \mathbf{0}_{n\times\left(n^{2}-n\right)}\\
\mathbf{0}_{\left(n^{2}-n\right)\times n} & \mathbf{0}_{\left(n^{2}-n\right)\times\left(n^{2}-n\right)}
\end{array}\right)\right\Vert ^{2}=\left\Vert \sum_{n\le k<n^{2}}\left(\mu_{k}\mathbf{u}_{k}\right)\cdot\left(\nu_{k}\mathbf{v}_{k}\right)^{\dagger}\right\Vert ^{2}+
\]
\[
\left\Vert \sum_{0\le k<n}\left(\mu_{k}\mathbf{u}_{k}\circ\left(\begin{array}{c}
\mathbf{0}_{n\times1}\\
\mathbf{1}_{\left(n^{2}-n\right)\times1}
\end{array}\right)\right)\cdot\left(\nu_{k}\mathbf{v}_{k}\circ\left(\begin{array}{c}
\mathbf{0}_{n\times1}\\
\mathbf{1}_{\left(n^{2}-n\right)\times1}
\end{array}\right)\right)^{\dagger}\right\Vert ^{2}
\]
which upper bounds the truncation error of the spectral approximation.
We note that the approximation error may be further reduced by the
use of iterative procedures based on the polynomial constraints which
generate the elimination ideals. The recursive approximation scheme
presented here allows us to build up an approximation for the spectral
decomposition of an $n\times n$ matrix starting from the spectral
decomposition of it's ${n \choose 2}$ matrix minors of size $2\times2$
all the way up to the $n$ matrix minors of size $\left(n-1\right)\times\left(n-1\right)$
from which we deduce the sought after approximation of the spectral
decomposition. The approximation algorithms described here straightforwardly
extend to hypermatrices of all orders.

\section{\label{sec:From-spectral-to}Unitary and Hermitian hypermatrices.}

We describe the spectral decomposition of even order Hermitian hypermatrices.
We also introduce here unitary hypermatrices and show how they can
be used to extend to hypermatrices the self adjoint argument for establishing
the existence of real solutions to spectral constraints. The discussion
here relates the spectral decomposition to the tensor eigenvalue first
defined by Lim \cite{Lim} and Qi \cite{Qi:2005:ERS:1740736.1740799}.
An even order hypermatrix $\mathbf{A}$ is said to be Hermitian if
$\overline{\mathbf{A}^{T}}=\mathbf{A}$. The general spectral decomposition
of fourth order hypermatrices is therefore expressed by 
\[
\mathbf{A}=\mbox{Prod}\left(\mbox{Prod}\left(\mathbf{Q},\boldsymbol{\Lambda}_{1},\boldsymbol{\Lambda}_{2},\boldsymbol{\Lambda}_{3}\right),\overline{\mbox{Prod}\left(\mathbf{U},\boldsymbol{\Gamma}_{1},\boldsymbol{\Gamma}_{2},\boldsymbol{\Gamma}_{3}\right)^{T^{3}}},\right.
\]
\[
\left.\mbox{Prod}\left(\mathbf{V},\boldsymbol{\Theta}_{1},\boldsymbol{\Theta}_{2},\boldsymbol{\Theta}_{3}\right)^{T^{2}},\overline{\mbox{Prod}\left(\mathbf{W},\boldsymbol{\Xi}_{1},\boldsymbol{\Xi}_{2},\boldsymbol{\Xi}_{3}\right)^{T}}\:\right)
\]
and
\[
\left[\mbox{Prod}\left(\mathbf{Q},\overline{\mathbf{U}^{T^{3}}},\mathbf{V}^{T^{2}},\overline{\mathbf{W}^{T}}\right)\right]_{i_{0},i_{1},i_{2},i_{3}}=\begin{cases}
\begin{array}{cc}
1 & \mbox{ if }i_{0}=i_{1}=i_{2}=i_{3}\\
0 & \mbox{otherwise}
\end{array}\end{cases}
\]
where the entries of the scaling hypermatrices are given by : 
\[
\left[\boldsymbol{\Lambda}_{1}\right]_{i_{0}i_{1}i_{2}i_{3}}:=\delta_{i_{1}i_{2}}\lambda_{i_{1}i_{3}},\left[\boldsymbol{\Lambda}_{2}\right]_{i_{0}i_{1}i_{2}i_{3}}:=\delta_{i_{1}i_{3}}\lambda_{i_{1}i_{0}},\left[\boldsymbol{\Lambda}_{3}\right]_{i_{0}i_{1}i_{2}i_{3}}:=\delta_{i_{1}i_{0}}\lambda_{i_{1}i_{2}}
\]
\[
\left[\boldsymbol{\Gamma}_{1}\right]_{i_{0}i_{1}i_{2}i_{3}}:=\delta_{i_{1}i_{2}}\gamma_{i_{1}i_{3}},\left[\boldsymbol{\Gamma}_{2}\right]_{i_{0}i_{1}i_{2}i_{3}}:=\delta_{i_{1}i_{3}}\gamma_{i_{1}i_{0}},\left[\boldsymbol{\Gamma}_{3}\right]_{i_{0}i_{1}i_{2}i_{3}}:=\delta_{i_{1}i_{0}}\gamma_{i_{1}i_{2}}
\]
\[
\left[\boldsymbol{\Theta}_{1}\right]_{i_{0}i_{1}i_{2}i_{3}}:=\delta_{i_{1}i_{2}}\theta_{i_{1}i_{3}},\left[\boldsymbol{\Theta}_{2}\right]_{i_{0}i_{1}i_{2}i_{3}}:=\delta_{i_{1}i_{3}}\theta_{i_{1}i_{0}},\left[\boldsymbol{\Theta}_{3}\right]_{i_{0}i_{1}i_{2}i_{3}}:=\delta_{i_{1}i_{0}}\theta_{i_{1}i_{2}}
\]
\[
\left[\boldsymbol{\Xi}_{1}\right]_{i_{0}i_{1}i_{2}i_{3}}:=\delta_{i_{1}i_{2}}\xi_{i_{1}i_{3}},\left[\boldsymbol{\Xi}_{1}\right]_{i_{0}i_{1}i_{2}i_{3}}:=\delta_{i_{1}i_{3}}\xi_{i_{1}i_{0}},\left[\boldsymbol{\Xi}_{3}\right]_{i_{0}i_{1}i_{2}i_{3}}:=\delta_{i_{1}i_{0}}\xi_{i_{1}i_{2}}
\]
where
\[
\delta_{i,j}=\begin{cases}
\begin{array}{cc}
1 & \mbox{if }0\le i=j<n\\
0 & \mbox{otherwise}
\end{array}\end{cases}.
\]
Entry wise the constraints are expressed as 
\[
a_{i_{0}i_{1}i_{2}i_{3}}=\left\langle \left(\boldsymbol{\mathbf{\lambda}}_{i_{0}}\circ\boldsymbol{\mathbf{\lambda}}_{i_{2}}\circ\boldsymbol{\mathbf{\lambda}}_{i_{3}}\right)\circ\mathbf{q}_{i_{0},i_{2},i_{3}},\overline{\left(\boldsymbol{\mathbf{\gamma}}_{i_{1}}\circ\boldsymbol{\mathbf{\gamma}}_{i_{0}}\circ\boldsymbol{\mathbf{\gamma}}_{i_{3}}\right)\circ\mathbf{u}_{i_{1}i_{3}i_{0}}},\right.
\]
\[
\left.\left(\boldsymbol{\theta}_{i_{2}}\circ\boldsymbol{\theta}_{i_{0}}\circ\boldsymbol{\theta}_{i_{1}}\right)\circ\mathbf{v}_{i_{2}i_{0}i_{1}},\,\overline{\left(\boldsymbol{\xi}_{i_{3}}\circ\boldsymbol{\xi}_{i_{1}}\circ\boldsymbol{\xi}_{i_{2}}\right)\circ\mathbf{w}_{i_{3}i_{1}i_{2}}}\:\right\rangle 
\]
\[
\delta_{i_{0}i_{1}i_{2}i_{3}}=\left\langle \mathbf{q}_{i_{0}i_{2}i_{3},},\overline{\mathbf{u}_{i_{1}i_{3}i_{0}}},\mathbf{v}_{i_{2}i_{0}i_{1}},\overline{\mathbf{w}_{i_{3}i_{1}i_{2}}}\right\rangle .
\]
In particular for unitary decomposition we have 
\[
a_{i_{0}i_{1}i_{2}i_{3}}=\left\langle \left(\boldsymbol{\mathbf{\lambda}}_{i_{0}}\circ\boldsymbol{\mathbf{\lambda}}_{i_{2}}\circ\boldsymbol{\mathbf{\lambda}}_{i_{3}}\right)\circ\mathbf{q}_{i_{0}i_{2}i_{3}},\overline{\left(\boldsymbol{\mathbf{\gamma}}_{i_{1}}\circ\boldsymbol{\mathbf{\gamma}}_{i_{0}}\circ\boldsymbol{\mathbf{\gamma}}_{i_{3}}\right)\circ\mathbf{q}_{i_{1}i_{3}i_{0}}},\right.
\]
\[
\left.\left(\boldsymbol{\theta}_{i_{2}}\circ\boldsymbol{\theta}_{i_{0}}\circ\boldsymbol{\theta}_{i_{1}}\right)\circ\mathbf{q}_{i_{2}i_{0}i_{1}},\overline{\left(\boldsymbol{\xi}_{i_{3}}\circ\boldsymbol{\xi}_{i_{1}}\circ\boldsymbol{\xi}_{i_{2}}\right)\circ\mathbf{q}_{i_{3}i_{1}i_{2}}}\:\right\rangle 
\]
and the unitarity constraints are entry wise expressed by 
\[
\delta_{i_{0},i_{1},i_{2},i_{3}}=\left\langle \mathbf{q}_{i_{0}i_{2}i_{3},},\overline{\mathbf{q}_{i_{1}i_{3}i_{0}}},\mathbf{q}_{i_{2}i_{0}i_{1}},\overline{\mathbf{q}_{i_{3}i_{1}i_{2}}}\right\rangle .
\]
The vectors $\left\{ \boldsymbol{\mathbf{\lambda}}_{i}\circ\boldsymbol{\mathbf{\lambda}}_{j}\circ\boldsymbol{\mathbf{\lambda}}_{k},\,\overline{\boldsymbol{\mathbf{\gamma}}_{i}\circ\boldsymbol{\mathbf{\gamma}}_{j}\circ\boldsymbol{\mathbf{\gamma}}_{k}},\,\boldsymbol{\theta}_{i}\circ\boldsymbol{\theta}_{i_{0}}\circ\boldsymbol{\theta}_{i_{1}},\overline{\boldsymbol{\xi}_{i_{3}}\circ\boldsymbol{\xi}_{i_{1}}\circ\boldsymbol{\xi}_{i_{2}}}\right\} _{0\le k<n}$
denote the vectors collecting the scaling values of the hypermatrix
$\mathbf{A}$. Note that the unitary decomposition described here
for even order hypermatrices is analogous to spectral decomposition
of Hermitian matrices expressed by 
\[
\mathbf{A}=\left(\mathbf{U}\cdot\mbox{diag}\left\{ \boldsymbol{\mu}\right\} \right)\cdot\overline{\left(\mathbf{U}\cdot\mbox{diag}\left\{ \boldsymbol{\nu}\right\} \right)^{T}},\quad\mathbf{U}\cdot\overline{\mathbf{U}^{T}}=\mathbf{I}_{n}
\]
entry wise expressed as 
\[
a_{i_{0}i_{1}}=\left\langle \boldsymbol{\mu}\circ\mathbf{u}_{i_{0}},\overline{\boldsymbol{\nu}\circ\mathbf{u}_{i_{1}}}\right\rangle ,\quad\delta_{i_{0}i_{1}}=\left\langle \mathbf{u}_{i_{0}},\overline{\mathbf{u}_{i_{1}}}\right\rangle 
\]
where the vectors $\boldsymbol{\mu},\overline{\boldsymbol{\nu}}$
correspond to the scaling vectors. Moreover, $\mathbf{A}$ is said
to admit slice invariant unitary decomposition if 
\[
\forall\:0\le i<j<n,\quad\boldsymbol{\mathbf{\lambda}}_{i}=\boldsymbol{\mathbf{\lambda}}_{j};\;\boldsymbol{\mathbf{\gamma}}_{i}=\boldsymbol{\mathbf{\gamma}}_{j};\;\boldsymbol{\theta}_{i}=\boldsymbol{\theta}_{i};\;\boldsymbol{\xi}_{i}=\boldsymbol{\xi}_{j}.
\]
In the case of matrices the spectral decomposition of a Hermitian
matrix is always slice invariant because the scaling vectors do not
change as the index of the eigenvector entries changes.
\begin{thm}
\label{thm:Self_adjointness}Let $\mathbf{A}$ denotes a Hermitian
hypermatrix which admits a slice invariant unitary decomposition then
it follows that the Hadamard product of the scaling vectors must be
real.
\end{thm}
The general argument of the proof is well illustrated for hypermatrices
of order $2$ and $4$. It will immediately be apparent how to extend
the argument to arbitrary even order hypermartices.
\begin{proof}
In the case of matrices we consider the bilinear form $\left\langle \mathbf{x},\,\overline{\mathbf{y}}\right\rangle _{\mathbf{A}}$
associated with the matrix $\mathbf{A}$. Let the spectral decomposition
of the matrix $\mathbf{A}$ be 
\[
\mathbf{A}=\left(\mathbf{U}\cdot\mbox{diag}\left\{ \boldsymbol{\mu}\right\} \right)\cdot\overline{\left(\mathbf{U}\cdot\mbox{diag}\left\{ \boldsymbol{\nu}\right\} \right)}^{T},\quad\mathbf{U}\cdot\overline{\mathbf{U}^{T}}=\mathbf{I}_{n}
\]
then the corresponding bilinear form can be expressed as follows 
\[
\left\langle \mathbf{x},\,\overline{\mathbf{y}}\right\rangle _{\mathbf{A}}=\sum_{0\le k<n}\left\langle \boldsymbol{\mu}\circ\mathbf{x},\,\overline{\boldsymbol{\nu}}\circ\overline{\mathbf{y}}\right\rangle _{\mathbf{u}_{k}\overline{\mathbf{u}_{k}^{T}}}
\]
where $\mathbf{u}_{k}$ denotes the $k$-th column of the unitary
matrix $\mathbf{Q}$. The bilinear form associated with the matrix
$\overline{\mathbf{A}^{T}}$ is therefore given by 
\[
\left\langle \mathbf{x},\,\overline{\mathbf{y}}\right\rangle _{\overline{\mathbf{A}^{T}}}=\sum_{0\le k<n}\left\langle \overline{\boldsymbol{\mu}}\circ\mathbf{x},\,\boldsymbol{\nu}\circ\overline{\mathbf{y}}\right\rangle _{\mathbf{u}_{k}\cdot\overline{\mathbf{u}_{k}^{T}}}.
\]
By Hermicity of $\mathbf{A}$ we have 
\[
\forall\mathbf{x},\mathbf{y}\in\mathbb{C}^{n},\quad\left\langle \mathbf{x},\,\overline{\mathbf{y}}\right\rangle _{\mathbf{A}}=\left\langle \mathbf{x},\,\overline{\mathbf{y}}\right\rangle _{\overline{\mathbf{A}^{T}}}\Rightarrow\boldsymbol{\mu}\circ\overline{\boldsymbol{\nu}}=\overline{\boldsymbol{\mu}}\circ\boldsymbol{\nu},
\]
by the combinatorial Nullstellensatz argument thus deriving that the
eigenvalues of $\mathbf{A}$ must all be real. Similarly we consider
the multilinear form associated with Hermitian hypermatrix $\mathbf{A}$
which admits a scale invariant unitary decomposition. The corresponding
multilinear form is expressed by 
\[
\left\langle \mathbf{x},\overline{\mathbf{y}},\mathbf{z},\overline{\mathbf{t}}\right\rangle _{\mathbf{A}}=\sum_{0\le k<n}\left\langle \left(\boldsymbol{\lambda}\circ\mathbf{x}\right),\left(\overline{\boldsymbol{\gamma}}\circ\overline{\mathbf{y}}\right),\left(\boldsymbol{\theta}\circ\mathbf{z}\right),\left(\overline{\boldsymbol{\xi}}\circ\overline{\mathbf{t}}\right)\right\rangle _{\mbox{Prod}\left(\mathbf{U}_{k},\overline{\mathbf{U}}_{k}^{T^{3}},\mathbf{U}_{k}^{T^{2}},\overline{\mathbf{U}}_{k}^{T}\right)}.
\]
The multilinear form associated with the hypermatrix $\overline{\mathbf{A}^{T}}$
is therefore given by 
\[
\left\langle \mathbf{x},\overline{\mathbf{y}},\mathbf{z},\overline{\mathbf{t}}\right\rangle _{\overline{\mathbf{A}^{T}}}=\sum_{0\le k<n}\left\langle \left(\boldsymbol{\gamma}\circ\mathbf{x}\right),\left(\overline{\boldsymbol{\theta}}\circ\overline{\mathbf{y}}\right),\left(\boldsymbol{\xi}\circ\mathbf{z}\right),\left(\overline{\boldsymbol{\lambda}}\circ\overline{\mathbf{t}}\right)\right\rangle _{\mbox{Prod}\left(\mathbf{U}_{k},\overline{\mathbf{U}}_{k}^{T^{3}},\mathbf{U}_{k}^{T^{2}},\overline{\mathbf{U}}_{k}^{T}\right)}.
\]
By Hermicity we have 
\[
\forall\,\mathbf{x},\mathbf{y},\mathbf{z},\mathbf{t}\in\mathbb{C}^{n},\quad\left\langle \mathbf{x},\overline{\mathbf{y}},\mathbf{z},\overline{\mathbf{t}}\right\rangle _{\mathbf{A}}=\left\langle \mathbf{x},\overline{\mathbf{y}},\mathbf{z},\overline{\mathbf{t}}\right\rangle _{\overline{\mathbf{A}^{T}}}\Rightarrow\boldsymbol{\lambda}\circ\overline{\boldsymbol{\gamma}}\circ\boldsymbol{\theta}\circ\overline{\boldsymbol{\xi}}=\overline{\boldsymbol{\lambda}}\circ\boldsymbol{\gamma}\circ\overline{\boldsymbol{\theta}}\circ\boldsymbol{\xi}.
\]

\end{proof}
\prettyref{thm:Self_adjointness} which extends to hypermatrices the
self-adjointness argument for establishing the existence of real solutions
to spectral constraints. We may also write that 
\[
\left\langle \mathbf{x},\overline{\mathbf{x}},\mathbf{x},\overline{\mathbf{x}}\right\rangle _{\mathbf{A}}=\sum_{0\le k<n}\left\langle \left(\boldsymbol{\lambda}\circ\mathbf{x}\right),\left(\overline{\boldsymbol{\gamma}}\circ\overline{\mathbf{x}}\right),\left(\boldsymbol{\theta}\circ\mathbf{x}\right),\left(\overline{\boldsymbol{\xi}}\circ\overline{\mathbf{x}}\right)\right\rangle _{\mbox{Prod}\left(\mathbf{U}_{k},\overline{\mathbf{U}}_{k}^{T^{3}},\mathbf{U}_{k}^{T^{2}},\overline{\mathbf{U}}_{k}^{T}\right)}.
\]
and therefore if we further make the simplifying assumption that for
some positive real number $\mu$ 
\[
\forall\,0\le i,\,k<n,\quad\mu\le\min\left\{ \left\langle \mathbf{e}_{i},\boldsymbol{\lambda}\right\rangle ,\,\left\langle \mathbf{e}_{i},\boldsymbol{\gamma}\right\rangle ,\,\left\langle \mathbf{e}_{i},\boldsymbol{\theta}\right\rangle ,\,\left\langle \mathbf{e}_{i},\boldsymbol{\xi}\right\rangle \right\} 
\]
and
\[
\max\left\{ \left\langle \mathbf{e}_{i},\boldsymbol{\lambda}\right\rangle ,\,\left\langle \mathbf{e}_{i},\boldsymbol{\gamma}\right\rangle ,\,\left\langle \mathbf{e}_{i},\boldsymbol{\theta}\right\rangle ,\,\left\langle \mathbf{e}_{i},\boldsymbol{\xi}\right\rangle \right\} \le\nu
\]
the entries of scaling hypermatrix therefore yield upper and lower
bounds for the eigenvalues for the symmetrized hypermatrix associated
with the multilinear forms $\left\langle \mathbf{x},\overline{\mathbf{x}},\mathbf{x},\overline{\mathbf{x}}\right\rangle _{\mathbf{A}}$
introduced by \cite{Lim,Qi:2005:ERS:1740736.1740799}. The corresponding
bounds are expressed by the inequality
\begin{equation}
\left\Vert \mu\,\mathbf{x}\right\Vert _{\ell_{4}}^{4}\le\left\langle \mathbf{x},\overline{\mathbf{x}},\mathbf{x},\overline{\mathbf{x}}\right\rangle _{\mathbf{A}}\le\left\Vert \nu\,\mathbf{x}\right\Vert _{\ell_{4}}^{4}.
\end{equation}
\bibliographystyle{amsalpha}
\bibliography{mybib}

\bigskip{}

\parbox[t]{1\columnwidth}{%
\noun{Department of Mathematics, Purdue University}\\
\noun{150 N. University Street, West Lafayette, IN 47907-2067}\\
E-mail:\texttt{ egnang@math.purdue.edu}%
}
\end{document}